\documentclass{amsart}
\usepackage{graphicx}
\newtheorem{theorem}{Theorem}
\newtheorem{lemma}[theorem]{Lemma}
\theoremstyle{definition}
\newtheorem{definition}[theorem]{Definition}
\theoremstyle{remark}
\newtheorem{remark}[theorem]{Remark}
\numberwithin{equation}{section}
\theoremstyle{corollary}

\theoremstyle{proposition}
\newtheorem{proposition}[theorem]{Proposition}

\newfont{\EUL}{eufm10 scaled 1000}

\newcommand\B{\mathbb{B}}
\newcommand\N{\mathbb{N}}

\newcommand\C{\mathbb{C}}

\newcommand\id{{\rm id}}

\begin{document} 
%
%
\title{Dynamics in the complex bidisc}
\author[Chiara Frosini]{Chiara Frosini$^\dag$}
\thanks{\rm $^\dag$ Supported by Progetto MURST di
Rilevante Interesse Nazionale {\it Propriet\`a geometriche delle
variet\`a reali e complesse} and by G.N.S.A.G.A (gruppo I.N.D.A.M)} 
\address{Dipartimento di
Matematica "U. Dini'', Universit\`a di Firenze, Viale Morgagni 67/A,
50134 Firenze , Italy.}
\email{frosini@math.unifi.it} \subjclass{Primary 32A40, 32H50.}
\date{January 15, 2004}
\keywords{Holomorphic maps, boundary behavior. }
 
\begin{abstract} 
Let  $\Delta^{n}$  be the unit polydisc in $\C^{n}$ and let $f$
be a holomorphic self map of $\Delta^{n}$.
When $n=1,$ it is well known, by Schwarz's
lemma, that $f$ has at most one fixed point in the unit disc. If no such point exists
 then $f$ has a unique boundary point, call it
$x\in\partial\Delta$, such that every horocycle $E(x,R)$ of center
$x$ and radius $R>0$ is sent into itself by $f$. This boundary point is
called the \emph{Wolff point of $f.$} In this paper we propose a  definition of Wolff points for holomorphic maps defined on a bounded domain of $\C^{n}$. In particular we characterize the set of Wolff points, $W(f),$ of a holomorphic self-map
$f$ of the bidisc in terms of the properties of the components of
the map $f$ itself.   
\end{abstract}
\maketitle

%
%
%
%
\section{Introduction}
Let  $D$ be  a bounded domain of $\C^{n}$ and let $f$
be a holomorphic self map of $D$.
We denote by $k_{D}$ the Kobayashi distance on $D$ ~\cite{Kob1},~\cite{Kob2},\cite{FranVes} and, as in \cite{AB1}, ~\cite{AB3}, we  define the \emph{small horosphere $E(x,R)$} and the \emph{big horosphere $F(x,R)$} of center $x$ and radius $R$ as follows:
\begin{equation}
\begin{array}{c}\label{eqHOR}                                                  E(x,R)=\{z\in D:\limsup\limits_{w\to x} [k_{D}(z,w)-k_{D}(0,w)]<\frac{1}{2}\log R\},\\
 F(x,R)=\{z\in D:\liminf\limits_{w\to x} [k_{D}(z,w)-k_{D}(0,w)]<\frac{1}{2}\log R\}.
\end{array}
\end{equation}
We say that
$\tau\in\partial D$ is a \emph{ Wolff point of f} if
$f(E(\tau,R))\subset E(\tau,R),$ for all $ R>0.$
Denote by $W(f)$ the set of \emph{Wolff points} of $f$
and denote by $T(f)$ the
\emph{ target set } 
\begin{center}$T(f):= \{ x \in \bar{D} \; | \exists \; \{ k_{n}\in \N
  \}, z \in \; D \; \hbox{such that} \; f^{k_{n}}(z)\to x \; \hbox{as}
  \; n\to \infty\}.$
\end{center}
If $D=\Delta=\{z\in\C:|z|<1\},$ it is well known, by Schwarz's lemma, that a holomorphic map $f:\Delta\to\Delta$
has at most one fixed point in $\Delta$. If $f$ has a fixed point in
$\Delta,$ say $z_{0},$ then $W(f)=\emptyset$ and $T(f)=\{z_{0}\}$. If $f$ has no fixed points in $\Delta$ then, by Wolff's lemma ~\cite{Wolff1}~\cite{RUDIN},
$W(f)$ is reduced to a boundary point $x\in \partial\Delta$ and, by the classical Denjoy theorem, $T(f)=W(f).$  
The same holds for self-maps of $\B^{n}$ with no interior fixed points ~\cite{HE1},~\cite{HE2},~\cite{MC1} and in particular if $f$ is a holomorphic self-map of a strongly convex domain
$D$ with $C^{3}-$boundary ~\cite{AB2}. 
On the other hand, if the map has fixed points in $D$ then either $W(f)=\emptyset$ and $T(f)$ is a unique point in $D,$ or $T(f)$ is a complex open subvariety $\Gamma$ of $D$ (affine in case $D=\B^{n}$) and $W(f)=\partial D\cup \bar{\Gamma}.$
In this paper we examine another type of convex domains (not strongly convex, not even with regular boundary): the
polydiscs. To avoid technical complications we restrict ourselves to dimension two, so we begin studying the case of the bidisc. 
In  this  case, we characterize
Wolff points of a holomorphic map $f:\Delta^{2}\to\Delta^{2}$  in terms
of the properties of the components of the map $f$.
Thus let $f:\Delta^{2}\to\Delta^{2}$ be a holomorphic
self-map in the complex bidisc without fixed points in
$\Delta^{2}$. Then $ f(x,y)=(f_{1}(x,y),f_{2}(x,y))$ with
$f_{1},f_{2}:\Delta^{2}\rightarrow \Delta$ holomorphic functions
in  $x$ and $y$. Then one of the two following possibilities holds (~\cite{HE1}):
\begin{enumerate}

      \item  there exists Wolff point of 
             $f_{1}(\cdot\;,y)$,  $e^{\imath\theta_{1}}$, which does  not depend on $y$
                                  or
      \item  there exists a holomorphic function
             $F_{1}:\Delta\rightarrow\Delta$, such that 
             $f_{1}(F_{1}(y),y)=F_{1}(y).$
             In this case  $f_{1}(x,y)= x \Rightarrow x=F_{1}(y).$\\
\end{enumerate}
Note that if $f\neq id_{\Delta}$ then cases $i),ii)$ exclude each
other.
Motivated by this result we  make the following definition:
\begin{definition}\label{MapType}
Let $f:\Delta^{2}\to\Delta^{2}$ be  a holomorphic function  and let
$f_{1},f_{2}$ be its components.
The map $f$ is called of:
\begin{enumerate}
\item \emph{first type} if:
      \begin{itemize}
       \item[-] there exists a holomorphic function
         $F_{1}:\Delta\rightarrow\Delta$, such that 
         $f_{1}(F_{1}(y),y)=F_{1}(y)$ and
       \item[-] there exists a holomorphic function
         $F_{2}:\Delta\rightarrow\Delta$, such that
          $f_{2}(x,F_{2}(x))=F_{2}(x).$
      \end{itemize}
\item \emph{second type} if (up to interchange $f_{1}$ with $f_{2}$):
      \begin{itemize}
       \item[-] there exists a  Wolff point of 
                $f_{1}(\cdot\;,y)$, $e^{\imath\theta_{1}}$, (necessarily independent of $y$) and
       \item[-] there exists a holomorphic function
                $F_{2}:\Delta\rightarrow\Delta$, such that 
                $f_{2}(x,F_{2}(x))=F_{2}(x).$
      \end{itemize}
\item \emph{third type} if:
      \begin{itemize}
       \item[-] there exists a Wolff point of 
                $f_{1}(\cdot\;,y)$, $e^{\imath\theta_{1}}$, (independent of $y$) and
       \item[-] there exists a Wolff point of
                $f_{2}(x,\cdot)$, $e^{\imath\theta_{2}}$, (independent of $x$).
      \end{itemize}
\end{enumerate}
\end{definition}
In case $f$ is of \emph{first type} and without interior fixed points
$F_{1}\circ F_{2}$ (respectively $F_{2}\circ F_{1}$) must have Wolff
point (see Lemma ~\ref{WolffF1F2}). Let  $e^{\imath\theta_{1}} $
(respectively by $e^{\imath\theta_{2}}$) be the Wolff point of  $F_{1}\circ F_{2}$ (respectively $F_{2}\circ F_{1}$.)
Let $\lambda_{1}:=\lim_{y\to e^{\imath\theta_{2}}}F'_{1}(y)$ and
$\lambda_{2}:=\lim_{x\to e^{\imath\theta_{1}}}F'_{2}(x)$ be
respectively  the \emph{boundary dilatation coefficients} of $F_{1}$
at $ e^{\imath\theta_{2}}$ and of $F_{2}$ at $ e^{\imath\theta_{1}}$
(see Lemma~\ref{WolffF1F2}). If $f$ is of second type we denote by
$\lambda_{2}$ the boundary dilatation coefficient of $F_{2}$ at
$e^{\imath\theta_{1}}.$ Finally  we let $\pi_{j}$ be the projection on
the $j-$component. With this notation our main result is:
\begin{theorem}\label{W(f)}
Let $f=(f_{1}, f_{2}) $ be a holomorphic map, without fixed
points in the complex bidisc. 
If $f_{1}\neq\pi_{1}$ and $f_{2}\neq\pi_{2},$  
then there are the following six cases:

$i)$ $\;\;W(f)=\emptyset$ if and only if $f$ is of \emph{first type} and
  $\lambda_{i}>1$ for one  $i=1,2.$

$ii)$ $\;W(f)={(e^{\imath{\theta_{1}}},e^{\imath{\theta_{2}}})}$ iff $f$ is of \emph{first type}
  $\lambda_{i}\leq 1$ for each  $i=1,2.$

$iii)$ $W(f)=\{\{e^{\imath{\theta_{1}}}\} \times \Delta\}\cup \{({(e^{\imath{\theta_{1}}},e^{\imath{\theta_{2}}})})\}$ iff $f$ is of
  \emph{second type} and
  $\lambda_{2}\leq 1.$

$iv)$ $\;W(f)=\{\{e^{\imath{\theta_{1}}}\} \times \Delta\}$ iff $f$ is of \emph{second type} and
  $\lambda_{2}> 1$.

$v)$ $\;\;W(f)= \{\{e^{\imath{\theta_{1}}}\}\times \Delta\}\cup \{{(e^{\imath{\theta_{1}}},e^{\imath{\theta_{2}}}})\}\cup \{\Delta \times
  \{e^{\imath{\theta_{2}}}\}\}$ iff $f$ is of \emph{third type}.\\
Otherwise if $f_{1}(x,y)=x$ $\forall\; y\in \Delta $ (or respectively
$f_{2}(x,y)=y$ $\forall\; x\in \Delta $ ) then:

$vi)$ $W(f)=\Delta^{2}\;\backslash \{\Delta\times\{e^{-\imath\theta_{2}}\}\}$
(or respectively $W(f)= \Delta^{2}\;\backslash \{\{e^{-\imath\theta_{1}}\}\times\Delta\}. $)
\end{theorem}
Let $f=(f_{1},f_{2}):\Delta^{2}\to\Delta^{2}$ be  holomorphic and with fixed points. 
By a result of Vigu\`e (see Proposition 4.1 ~\cite{VI1}) it follows that  $f_{1}( \cdot , y)$ and $f_{2}(x, \cdot)$ also must have fixed points. 
Notice that if $\dim Fix(f)=2$ then $f=\id|_{\Delta^{2}}.$ Moreover if $\dim Fix(f)=1$ we know (see Propositions 2.6.10; 2.6.24 ~\cite{AB2}) that $Fix (f)$ is a geodesic of $\Delta^{2}$ and then 
$ Fix(f)$ can be parametrized as $ \Delta\ni z \to(g(z),z)$ with $
g\in Hol(\Delta,\Delta).$

\begin{theorem}\label{FIX1}Let $f=(f_{1},f_{2}):\Delta^{2}\to\Delta^{2}$
  be  a holomorphic  map, not an automorphism, with fixed points in
$\Delta^2$. Assume, up to automorphisms, that $f(0,0)=(0,0)$.
\begin{itemize}
\item If $\dim Fix(f)=0$ then $W(f)=\emptyset.$
\item If $\dim Fix(f)=1$ and we let $G:= Fix(f) $ then :
\end{itemize}
\begin{enumerate}
\item $g(z\in Aut(\Delta)\cup\{id\}$  iff
$W(f)=\partial G$ (and this is the case iff there exists a point
$(e^{i\theta},1)\in(\partial \Delta)^2$ which belongs to $W(f)$);
\item $g\not\in Aut(\Delta)\cup\{id\}$  is a proper map iff
$W(f)=\emptyset;$
\item $g$ is not a proper map iff $W(f)$ is
disconnected (and this is the case iff $f_{2}= \pi_{2}$).
\end{enumerate}
\end{theorem}
If $f\in Aut (\Delta^{2})$ has fixed points in $\Delta^{2},$ its components are elliptic automorphisms of $\Delta$ and $W(f)=\emptyset.$

In  section  2 we are going to introduce some useful tools  to prove our main theorem.
In particular we describe the property of the component $f_{1}$ and $f_{2}$ of the function $f:\Delta^{2}\to \Delta^{2}$
and we also study some property of the set of\emph{Wolff points} of $f.$
Using these results in  section 3  we prove  Theorem~\ref{W(f)} and
Theorem ~\ref{FIX1}. 
Finally in section 4 we give some examples.  

I want to sincerely thank Filippo Bracci
 for  his  continuous  assistance and professor Graziano Gentili for
 many useful conversation  for this
 work. Also, I wish to thank professor Pietro Poggi-Corradini for his
 helpful comments concerning this work.


%
%
%


\section{Preliminary results}
We need now to introduce some notation and some preliminary results,
as a generalization of  Julia's  lemma due to Abate (Abate~\cite{AB1}):
\begin{theorem}(Abate~\cite{AB1})~\label{JuliaForPolydiscs}
Let  $f:\Delta^{n}\to\Delta $ be a holomorphic map and let  $x\in
\partial\Delta^{n}$ be such that
\begin{equation}~\label{limiteJulia}
\liminf_{\tilde{w}\to
  x}[k_{\Delta^{n}}(0,\tilde{w})-\omega(0,f(\tilde{w}))]=
\frac{1}{2} \log\alpha_{f} <\infty
\end{equation}
then  there exists $\tau\in \partial\Delta $ such that $\forall \;\; R>0, $
$f(E(x,R))\subseteq E(\tau,\alpha_{f} R).$
Furthermore  $f$ admits restricted $E-limit $  $\tau $ at $x.$ (see
(Abate~\cite{AB1}) for the definition).
\end{theorem}
Consider $f=(f_{1}, f_{2})\in Hol(\Delta^{2}, \Delta^{2})$.
Using the result of Herv{\'e} [see Theorem 1 in  ~\cite{HE1}] and  Definition ~\ref{MapType}, we are going to examine the properties of $f,\; f_{1}\;\hbox{and}\; f_{2}$. 
\begin{lemma}~\label{WolffF1F2}
  Let $f=(f_{1}, f_{2})\in Hol(\Delta^{2}, \Delta^{2}).$ If $f$ is a map of \emph{first type} then:
\begin{itemize}
\item[i)] The function 
$F_{1}\circ F_{2}$ (respectively $F_{2}\circ F_{1}$) has Wolff point $e^{\imath\theta_{1}}$ (respectively $e^{\imath\theta_{2}}$). In this case we let  $\lambda_{12}$ (respectively $\lambda_{21}$) be  the boundary dilatation coefficient of $F_{1}\circ F_{2}$ (resp. $F_{2}\circ F_{1}$) at its Wolff point;
\item[ii)] $F_{1} $ has non-tangential limit $e^{\imath\theta_{1}}$ at
      $e^{\imath\theta_{2}}$ and
      \begin{equation} \liminf_{y\to
      e^{\imath\theta_{2}}}  \frac{1-|F_{1}(y)|}{1-|y|}=\lambda_{1}
    \;\;\; ;\;
      0< \lambda_{1} < +\infty ;\end{equation}
\item[iii)] $F_{2} $ has non-tangential limit $e^{\imath\theta_{2}}$ at
      $e^{\imath\theta_{1}}$ and
      \begin{equation} \liminf_{x\to
      e^{\imath\theta_{1}}}  \frac{1-|F_{2}(x)|}{1-|x|}=\lambda_{2}
    \;\;\; ;\; 0< \lambda_{2} < +\infty. \end{equation}
\end{itemize}
Furthermore $\lambda_{12}=\lambda_{21}=\lambda_{1}\lambda_{2}.$
\end{lemma}
\begin{proof}
$i).$
We note that if $f$ is a map of \emph{first type} then the function
$F_{1}\circ F_{2}$ cannot have fixed points in $\Delta$ or otherwise we can build fixed points in $\Delta^{2}$ for $f.$ Infact,
suppose that  $x_{0}\in \Delta$ is such that  $(F_{1}\circ
F_{2})(x_{0})\nolinebreak =\nolinebreak x_{0}$ and let $F_{2}(x_{0})=y_{0}.$ Then 
$f(x_{0},y_{0})=(f_{1}(F_{1}(F_{2}(x_{0})),F_{2}(x_{0})),f_{2}(x_{0},F_{2}(x_{0})))=(x_{0},y_{0}).$
By   Wolff's lemma, $F_{1}\circ F_{2}$ has a Wolff point, say $e^{\imath\theta_{1}}\in\partial\Delta.$
Similarly $F_{2}\circ F_{1}$ has a Wolff
point, say $e^{\imath\theta_{2}}\in\partial\Delta.$
It is well known, by  classical results of  Julia, Wolff and Caratheodory ~\cite{AB2}, that:
 $\liminf_{x\to
           e^{\imath\theta_{1}}}\frac{1-|F_{1}(F_{2}(x))|}
           {1-|x|}=\lambda_{12} \;\;
          \hbox{for some}\;  0<\lambda_{12}\leq 1 ;$ and         
 $ \liminf_{y\to
          e^{\imath\theta_{2}}}
        \frac{1-|F_{2}(F_{1}(y))|}{1-|y|}=\lambda_{21}\;\;
     \hbox{for some}\;     0<\lambda_{21}\leq 1.$
         
$ii)$
By Schwarz's lemma we know that 
$\frac{1-|F_{1}(y)|}{1-|y|}\geq \frac{1-|F_{1}(0)|}{1+|F_{1}(0)|}:=M$ for $y\in\Delta.$
If we take $y=F_{2}(x),$ it implies:
$ 1- |F_{1}(F_{2}(x))|\geq M(1-|F_{2}(x)|)$
and dividing by $(1-|x|)$ we obtain:
$\frac{1-|F_{1}(F_{2}(x))|}{1-|x|} \geq M \frac{1-|F_{2}(x)|}{1-|x|}
\geq M \frac{1-|F_{2}(0)|}{1+|F_{2}(0)|}.$
Let   $ \lambda_{2}:= \liminf_{x\to e^{\imath\theta_{1}}}
\frac{1-|F_{2}(x)|}{1-|x|} $ be the boundary dilatation coefficient of $ F_{2}
$
in $ e^{\imath\theta_{1}}.$
Thus  we can conclude that:
$ 1\geq\lambda_{12}\geq M \lambda_{2} >0. $
Note that, since $ M\leq 1, $ it follows that  $\lambda_{2}$ is a
finite and positive number. But we also know  that $|F_{2}(x)|\to 1 $
for $x\to e^{\imath\theta_{1}}$ and by Julia's lemma we can
conclude that there exists a unique $ e^{\imath\gamma_{2}}$ such that the non-tangential  limit
 of $F_{2}$ in $e^{\imath\theta_{1}}$ is equal to $ e^{\imath\gamma_{2}}$, that is:
$ K-\lim_{x\to
e^{\imath\theta_{1}}} F_{2}(x)= e^{\imath\gamma_{2}}.$ Let us consider the function $ (F_{2}\circ F_{1}):\Delta\to\Delta$. As we
proved in $i)$, this  function has Wolff point  $e^{\imath\theta_{2}}.$ Notice that, as for $\lambda_{2},$ one can prove that $\lambda_{1}:=\liminf_{y\to e^{\imath\theta_{2}}}F_{1}(y)$ is such that $0<\lambda_{1}<+\infty.$
We are going to show that:
$e^{\imath\gamma_{2}}=e^{\imath\theta_{2}}.$
Let us consider the sequence $\{(F_{2}\circ F_{1})^{k}\}$ of iterates of
$(F_{2}\circ F_{1})$. By the Wolff-Denjoy theorem it converges, uniformly
on compact sets, to the constant  $e^{\imath\theta_{2}}.$
Notice that
$(F_{2}\circ F_{1})^{k}=F_{2}\circ (F_{1}\circ F_{2})^{k-1}\circ F_{1}$
and, furthermore, $e^{\imath\theta_{1}}$ is the Wolff point of $(F_{1}\circ
F_{2}).$ Therefore
$(F_{1}\circ F_{2})^{k}\to
e^{\imath\theta_{1}} $ as $k\to\infty.$
Let us fix  $y_{0}\in \Delta$ and let  $x_{0}=F_{1}(y_{0}).$
Set $w_{k}:=(F_{1}\circ F_{2})^{k}(x_{0})$ and $z_{k}:=(F_{2}\circ F_{1})^{k}(y_{0})$.
Then we have $ w_{k}\to e^{\imath\theta_{1}}$ and $z_{k}\to e^{\imath\theta_{2}}$ as $k\to \infty.$
Moreover we notice that:
$\omega(0, w_{k})-\omega(0,F_{2}(w_{k}))=
 \omega (0,F_{1}( z_{k}))-\omega(0,z_{k+1})= \omega (0,F_{1}( z_{k}))-\omega(0,z_{k})+\omega(0,z_{k})-\omega(0,z_{k+1}).$
Then
\begin{center}$
-\limsup_{k\to\infty}[\omega (0, F_{1}(z_{k}))(y_{0}))-
\omega (0,
 z_{k}(y_{0}))]=\liminf_{k\to\infty} -[\omega (0, F_{1}(z_{k}))(y_{0}))-\omega (0,
 z_{k}(y_{0}))]$
$
\geq\liminf_{w\to e^{\imath\theta_{2}}} -[\omega (0,
F_{1}(w))-\omega (0,w)]
=\frac{1}{2}\log \lambda_{1}>-\infty,
$
\end{center}
we conclude that
\begin{center}
$\limsup_{k\to\infty}\omega(0, (F_{1}\circ
  F_{2})^{k}(x_{0}))-\omega(0,F_{2}(F_{1}\circ
  F_{2})^{k}(x_{0}))$
$=\limsup_{k\to\infty}[\omega (0, F_{1}((F_{2}\circ F_{1})^{k})(y_{0}))-\omega (0,
  (F_{2}\circ F_{1})^{k}(y_{0}))]<+\infty.$
\end{center}
Furthermore, since $(F^{2}\circ F^{1})^{k}$  converges uniformly on
  compact sets to its Wolff point we have that:
\begin{center}$\limsup_{k\to\infty}[\omega (0,
 z_{k}(y_{0}))- \omega (0,
 z_{k+1}(y_{0}))]=\limsup_{k\to\infty}[\omega (0,
 (F_{2}\circ F_{1})^{k}(y_{0}))- \omega (0,
 (F_{2}\circ F_{1})^{k+1}(y_{0}))]$
$\leq \limsup_{k\to\infty}|\frac{1+|(F_{2}\circ F_{1})^{k}(y_{0})|}{1+(F_{2}\circ F_{1})^{k+1}(y_{0})}\frac{1-|(F_{2}\circ F_{1})^{k+1}(y_{0})|}{1-|(F_{2}\circ F_{1})^{k}(y_{0})|}|\leq 0.$
\end{center}
We have 
$\lim_{k\to\infty}F_{2}((F_{1}\circ
F_{2})^{k}(x_{0}))=e^{\imath\gamma_{2}}$ and 
$e^{\imath\gamma_{2}}=\lim_{k\to\infty}F_{2}((F_{1}\circ
F_{2})^{k}(x_{0}))\\=\lim_{k\to\infty}(F_{2}\circ
F_{1})^{k}(y_{0}))=e^{\imath\theta_{2}}$
thus  we can conclude that
$e^{\imath\gamma_{2}}=e^{\imath\theta_{2}}.$

$iii).$
The last thing we need to prove is that  $F_{1}$ has  non tangential limit
$e^{\imath\theta_{1}}$ at $e^{\imath\theta_{2}}.$ This can be done
proceeding as before.
Using Julia-Wolff-Caratheodory theorem~\cite{Wolff2},~\cite{AB2} we obtain that
$\lambda_{1}=\frac{\lambda_{12}}{\lambda_{2}},$ $\lambda_{2}=\frac{\lambda_{21}}{\lambda_{1}}$ 
and, in particular
$
\lambda_{21}=\lambda_{1} \lambda_{2}= \lambda_{12}. 
$
\end{proof}
It's clear that $\lambda_{1}
>1 $ implies $\lambda_{2}<1$ (and  the converse is also true.).
To simplify notations,
up to automorphisms, from now on assume
$e^{\imath\theta_{1}}=e^{\imath\theta_{2}}=1.$
Let  $\alpha_{f_{1}}=$ denote the number defined in theorem~\ref{JuliaForPolydiscs} for the function $f_{1}$ at the point $(1,1):$

\begin{lemma}\label{DilatCoeffHoros}
Let $f=(f_{1}, f_{2}):\Delta^{2}\to\Delta^{2}$ be  holomorphic and suppose that  $f_{1}\neq\pi_{1}$ and
$f_{2}\neq\pi_{2}$. 
If $f_{1}(E((1,1),R))\subseteq E(1, R)$ $\forall \; R>0$  then $\alpha_{f_{1}}\leq 1.$
\end{lemma}
\begin{proof}
Let  consider the holomorphic function  $\varphi:\Delta\to\Delta
$ defined by $\varphi(\xi)=f_{1}(\xi,\xi)$
and let $\alpha_{\varphi}$ be the boundary dilatation coefficient of the map $\varphi$ at the point $1.$
Let $|||(z,w)|||:=\max\{||z||, ||w||\}.$
Since (see lemma $3.2$ ~\cite{AB1})
$\liminf_{t\to
    1^{-}}\frac{1-|f_{1}(t,t)|}{1-|t|}=\alpha_{f_{1}}$
then we have that
$\alpha_{f_{1}}=\liminf_{(z,w)\to
    (1,1)}\frac{1-|f_{1}(z,w)|}{1-|||(z,w)|||}\leq \liminf_{\xi\to
    1}\frac{1-|f_{1}(\xi,\xi)|}{1-|\xi|}
\leq \liminf_{t\to
    1^{-}}\frac{1-|f_{1}(t,t)|}{1-|t|}=\nolinebreak[4]\alpha_{f_{1}}.$
And
$\alpha_{\varphi}=\liminf_{\xi\to 1}\frac{1-|\varphi(\xi)|}{1-|\xi|}=\liminf_{\xi\to 1}\frac{1-|f_{1}(\xi,\xi)|}{1-|\xi|}=\liminf_{t\to
    1^{-}}\frac{1-|f_{1}(t,t)|}{1-t}=\alpha_{f_{1}}$
and we can conclude that $\alpha_{\varphi}=\alpha_{f_{1}}.$
But since  $f_{1}(E((1,1),R))\subseteq E(1, R)$ $\forall \; R>0$ then $\alpha_{\varphi}\leq 1$ and this ends the proof.
\end{proof} 

\begin{proposition}~\label{R1R2}
Let $f=(f_{1},f_{2}):\Delta^{2} \to \Delta^{2} $ be  holomorphic and without
fixed points in $\Delta^{2}.$ Suppose that  $f_{1}\neq\pi_{1}$ and
$f_{2}\neq\pi_{2}$. Fix $R_{1},R_{2}>0$ such that 
$\frac{1}{R_{1}}=\frac{\lambda_{2}}{R_{2}}.$ If either  
\begin{enumerate}
\item[i)] $f$ is of \emph{first type} and $\lambda_{1}>1$ or
\item[ii)] $f$ is of \emph{second type} and $\lambda_{2}>1,$ 
\end{enumerate} 
Then
$f(E(1,R_{1})\times E(1,R_{2}))\subseteq E(1,R_{1})\times E(1,
R_{2}).$
\end{proposition}
\begin{proof}

i)
Let consider $(x_{0},y_{0})\in E(1,R_{1})\times E(1, R_{2}).$ 
We have:
\begin{center}$\lim\limits_{z\to 1}
\omega(f_{1}(x_{0},y_{0}),z)-\omega(0,z)= \lim\limits_{t\to
  1^{-}}\omega(f_{1}(x_{0},y_{0}),f_{1}(F_{1}(F_{2}(t)),
F_{2}(t)))-\omega(0,F_{1}(F_{2}(t)))$
$\leq \lim\limits_{t\to 1^{-}}\max\{\omega(x_{0},F_{1}(F_{2}(t))),
  \omega(y_{0},F_{2}(t))\}-\omega(0,F_{1}(F_{2}(t))).$
\end{center}
If $\max\{\omega(x_{0},F_{1}(F_{2}(t))),
  \omega(y_{0},F_{2}(t))\}=\omega(x_{0},F_{1}(F_{2}(t)))$
then 
$\lim_{t\to
  1^{-}}\omega(x_{0},F_{1}(F_{2}(t)))-\omega(0,F_{1}(F_{2}(t)))\leq
\frac{1}{2}\log R_{1}\;\hbox{because}\; x_{0}\in E(1,R_{1}).$
If
$\max\{\omega(x_{0},F_{1}(F_{2}(t))),
  \omega(y_{0},F_{2}(t))\}=\omega(y_{0},F_{2}(t))$
then:
\begin{center}
$\limsup\limits_{t\to
  1^{-}}\omega(y_{0},F_{2}(t))-\omega(0,F_{1}(F_{2}(t)))$\\
$\leq \limsup\limits_{t\to
  1^{-}}\omega(y_{0},F_{2}(t))-\omega(0,F_{2}(t))+\limsup\limits_{t\to
  1^{-}}\omega(0,F_{2}(t))-\omega(0,F_{1}(F_{2}(t))) $\\
$\leq\frac{1}{2}\log R_{2} +\limsup\limits_{t\to
  1^{-}}\frac{1}{2} \log [\frac{1+|F_{2}(t)|}{1+|F_{1}(F_{2}(t))|}\frac{|1-F_{1}(F_{2}(t))|}{|1-F_{2}(t)|}\frac{|1-F_{2}(t)|}{1-t}\frac{1-t}{1-|F_{2}(t)|}]$\\$\leq \frac{1}{2}\log R_{2}+ \frac{1}{2}\log
\lambda_{1}=\frac{1}{2}\log\lambda_{1}R_{2}=\frac{1}{2}\log\frac{\lambda_{12}}{\lambda{2}}R_{2}=\frac{1}{2}\log
R_{1},$ 
\end{center}
 by Julia-Wolff-Caratheodory theorem.
So $f_{1}(x_{0},y_{0})\in E(1, R_{1}).$
In the same way one can prove that $f_{2}(x_{0},y_{0})\in E(1, R_{2}).$
So $f(E(1,R_{1})\times E(1,R_{2}))\subseteq E(1,R_{1})\times E(1,R_{2}).$
This end the proof of the point \emph{i).}
\\
\emph{ii).}
By  hypothesis it  follows immediately that $R_{1}\geq R_{2}$ and then we have:
$$f_{1}(E(1,R_{1})\times (E(1,R_{2})))\subseteq f_{1}(E(1,R_{1})\times
\Delta)\subseteq E(1,R_{1}).$$
So it remains to examine the second component $f_{2}.$ Again by  hypothesis we
have that $\lambda_{2}=\frac{R_{2}}{R_{1}}$ and we can proceed as
above and we prove that  $f_{2}(E(1,R_{1})\times (E(1,R_{2})))\subseteq E(1,R_{2})$ and it concludes the proof of the point \nolinebreak[4]\emph{ii).} \end{proof}

\begin{proposition}~\label{W(f)Cpa}
Let $f=(f_{1}, f_{2}):\Delta^{2}\to\Delta^{2}$ be  holomorphic and without
fixed points in $\Delta^{2}.$ Suppose that  $f_{1}\neq\pi_{1}$ and
$f_{2}\neq\pi_{2}$. Then $W(f)$ is
arcwise connected.
\end{proposition}
\begin{proof}[Dim.]
 Let  $\Gamma_{(x,y)}$ denote the flat component of the
boundary of
$\Delta^{2}$ containing the  point $(x,y)\in\partial\Delta^{2}$. We can consider the
following cases:\\
$1)$ ``$z, w \in W(f)\cap \Gamma_{(0,1)}.$" 
In this case every point of the component $\Gamma_{(0,1)}$ is a Wolff
point. Indeed, if $x\in\Gamma_{(0,1)}$ then $F(x,R)=E(x,R)$ and in particular they coincide with $E(z,R)=F(z,R)=E(w,R)=F(w,R).$ Thus if a point of the flat component $\Gamma_{(0,1)}$ is a Wolff point then every point of this component is a Wolff point. We can conclude that there exists a continuous   path, of Wolff points, that links $z$ and $w$.\\
$2)$ ``$z\in W(f)\cap \partial\Gamma_{(0,1)}$ and $w\in
  W(f)\cap\Gamma_{(0,1)}$."
 $w\in W(f)$ implies that every point of the flat component  $\Gamma_{(0,1)}$
is a Wolff point then, it is sufficient consider
the radius linking $z$ and $w.$\\
$3)$ ``$ z\in W(f)\cap \Gamma_{(0,1)}$ and $w\in
  W(f)\cap\Gamma_{(1,0)}.$"
Every point of the flat components $\Gamma_{(0,1)}$ and $\Gamma_{(1,0)}$ is a Wolff point of $f$.
By definition of Wolff points and horosphere, it follows that:
$f(E(1,R) \times  E(1,R))= f((\Delta\times E(1,R))\cap (E(1,R)\times \Delta))
\subseteq E(1,R)\times E(1,R)
\;\forall\; R>0$
and  $(1,1)\in W(f)$; thus there is a continuous   path that links $z$ and $w,$ passing through the point $(1,1).$

We are now going to prove that there isn't any other possibility.
In fact the following cases are not possible:
\\
\emph{a)} ``$\exists\;z\in W(f)\cap\Gamma_{(0,1)}$ e $\exists w \in W(f)\cap
  \Gamma_{(0, -1)}.$"
In this case, by definition of Wolff point, we would have:
$f_{2}(\Delta\times E(1,R))\subseteq E(1,R) \; \forall\; R>0 \;\hbox{and}\;
f_{2}(\Delta\times E(-1,R))\subseteq E(-1,R) \; \forall\; R>0.
$
And if we consider  $ y_{0}\in E(1,R_{1})\cap E(-1,R_{2})$ and $x\in \Delta$ we
would  have:
$f_{2}(x,y_{0})=\nolinebreak[4]y_{0}$ $\forall\;\; x\in \Delta.$ But it is possible only if
$f_{2}(x,y)=y$ for all $ x\in \Delta$ and it isn't our case.
\\
\emph{b)} `` $\exists \; z\in W(f)\cap \partial\Gamma_{(0,1)} $ and $\exists \;
  w\in W(f)\cap \partial\Gamma_{(0,1)}$, $z\neq w.$"
We can suppose that $z=(-1,1)$ and $w=(1,1)$. Then
$f_{2}(\Delta\times E(1,R))\subseteq E(1,R) \; \forall\; R>0 \; \hbox{and}\;
f_{2}(\Delta\times E(-1,R))\subseteq E(-1,R) \; \forall\; R>0.
$
Again, if we take $y_{0}\in \partial E(1,R_{1}) \cap \partial E(-1,
R_{2}),$ chosen opportunely  $R_{1}, R_{2} >0,$ we obtain that $f_{2}(x, y_{0})= y_{0}$
$\forall\;\; x\in \Delta.$ But it is possible only if
$f_{2}(x,y)=y$ $\forall \; x\in \Delta$ and it isn't our case.
\\
\emph{c)} ``$(1,1) ; (-1,-1) \in W(f)"$.
In this case we have
$f(E(-1,R) \times E(-1,R))\subseteq  E(-1,R)\times E(-1,R)  \; \forall\; R>0 \; \hbox{and}\;
f(E(1,R) \times E(1,R))\subseteq E(1,R)\times E(1,R)  \; \forall\; R>0.$
Then if we take $0\in \partial E(1,R)\cap \partial E(-1,R)$ we obtain that 
$(0,0)\in [\partial(E(-1,R) \times E(-1,R))]\cap [\partial E(1,R)
\times E(1,R)]$ and consequently: $f(0,0)=0$, but it is not possible
because we supposed  $f$ without fixed points in the bidisc.
\end{proof}

%
%
%
%
\section{Wolff points}

We are finally ready to prove our main result , theorem ~\ref{W(f)}, taking
$e^{\imath\theta_{1}}=e^{\imath\theta_{2}}=1:$
\begin{proof}[Proof of Theorem ~\ref{W(f)}]

We are going to prove every statement in the direction ``
$\Leftarrow$.\nolinebreak[4]"

$ii).$
Let apply   Lemma ~\ref{JuliaForPolydiscs} to the holomorphic function
$f_{1}:\Delta^{2}\to\Delta$ and let us study:
\begin{equation}\label{Julia per f_{1}}
\liminf_{(x,y)\to(1,1)} K_{\Delta^{2}}((0,0),(x,y))-\omega(0,f_{1}(x,y))
\end{equation}
Consider the direction $x=F_{1}(w)\;\; y=w.$
We obtain:
$\omega(0,x)=\omega(0,F_{1}(w)) \;\hbox{and}\\\;\omega(0,y)=\omega(0,w).$
Since $\lambda_{1}\leq 1$, then:
$\liminf_{w\to 1} [\omega(0,w)-\omega(0,F_{1}(w))]\leq 0.$
In par\-ti\-cu\-lar there exists a subsequences $w_{k}$ such that $\lim_{k\to\infty}w_{k}= 1$ and
$\lim_{k\to
  +\infty}[\omega(0,w_{k})-\omega(0,F_{1}(w_{k}))]\leq 0.$
If we look at (~\ref{Julia per f_{1}}) along the direction $w_{k}$ we obtain:
\begin{center}$\liminf_{(x,y)\to
  (1,1)}K_{\Delta^{2}}((0,0),(x,y))-\omega(0,f_{1}(x,y))
\leq\liminf_{k\to
  +\infty}K_{\Delta^{2}}((0,0),(F_{1}(w_{k}),w_{k}))-\omega(0,f_{1}(F_{1}(w_{k}),w_{k}))=$\\
$=\liminf_{k\to
  +\infty}\omega(0,F_{1}(w_{k}))-\omega(0,F_{1}(w_{k}))=0 .$
\end{center}
Then, by (~\ref{Julia per f_{1}}), there exists $ \tau_{1}\in \partial\Delta$ such that
$f_{1}(E(1,1),R)\subseteq E(\tau_{1},R)$ $\forall \; R>0. $
In particular $f_{1}$ admits restricted E-limit  $\tau_{1}.$
Let apply again the ~\ref{Julia per f_{1}} to the map
$f_{2}:\Delta^{2}\to\Delta$ along the direction $x=z$ ; $y=F_{2}(z).$
Proceeding as above, we obtain
$\liminf_{(x,y)\to (1,1)}
K_{\Delta^{2}}((0,0),(x,y))-\omega(0,f_{2}(x,y))=0.$
Then there exists $ \tau_{2}\in
\partial\Delta $ such that $f_{2}(E(1,1),R)\subseteq E(\tau_{2},R)$
$\forall \; R>0 $ and $f_{2}$ admits restricted E-limit
$\tau_{2}.$
We claim that $\tau_{1}=\tau_{2}=1.$
Consider the curves $\sigma_{1}(t)=(F_{1}(t),t)$ and
$\sigma_{2}(t)=(t, F_{2}(t))$; $t\in [0,1). $
These curves, 
$\sigma_{1},\sigma_{2},$ are peculiar $(1,1)-$curves.
Furthermore
$ f_{1}(\sigma_{1}(t))=f_{1}(F_{1}(t),t)=F_{1}(t)\to\nolinebreak[4]1 \;\; \hbox{as}\;\;
t\to \nolinebreak[4]1^{-} $ and
$f_{2}(\sigma_{2}(t))=f_{2}(t, F_{2}(t))=F_{2}(t)\to 1 \;\; \hbox{as} \;\;
t\to 1^{-}.$
Since  $\tau_{1},\tau_{2},$ are respectively
restricted E-limit of $f_{1}$ and restricted E-limit of $f_{2}$
then we can conclude that
$\tau_{1}=1$ and $ f_{1}(E((1,1),R))\subseteq  E(1,R)\; \forall\; R>0$;
$\tau_{2}=1$ and $ f_{2}(E((1,1),R))\subseteq E(1,R) \; \forall\; R>0.$
Thus if $f$ is of \emph{first type} and $\lambda_{1}\leq 1$,
$\lambda_{2}\leq 1 $   then:
 $$f((E(1,1),R)=(f_{1}(E((1,1),R)), f_{2}(E((1,1),R)))\subseteq
 E((1,1),R)\;\;\forall\;R>0$$
and in particular \emph{$\tau={(1,1)}\in W(f).$}
Notice  that, by the proof of Lemma ~\ref{W(f)Cpa},
points as ${(-1,1)}$ and ${(1,-1)}$ cannot be Wolff point of $f$,
because they are on the Silov boundary of the same flat component
of the point ${(1,1)}$. Furthermore, also, no points of the flat
components $\Gamma_{(0,1)}$, $\Gamma_{(1,0)}$ can be  Wolff points of $f$
otherwise, by definition, the point $1$ would be a Wolff point for $f_{1}(\cdot,y)$ and $f_{2}(x, \cdot)$ and it isn't possible because
$f_{1}(\cdot, y)$ has fixed points described by the function
$F_{1}$ and $f_{2}(x, \cdot)$ has fixed points described by the
function $F_{2}$. By lemma  ~\ref{W(f)Cpa}, we know that $W(f)$ is
arcwise connected, then the unique Wolff point of $f$ must be
$(1,1)$.
\\
\emph{$ i)$}
Since $F_{1}$ is a holomorphic self-map of $\Delta,$ there are at most two possibilities  :
 $F_{1}(y_{0})=y_{0} $ for some $y_{0}\in \Delta$ or  $F_{1}$ has a Wolff point $\tau_{1}\neq 1.$ We can suppose,
  $\tau_{1}=-1.$
In this last case we have $K-\lim\limits_{y\to -1} F_{1}(y)=-1$ and the
angular derivative of $F_{1}$ in $\{-1\}$ is  $\delta_{1}\leq 1.$
Suppose $F_{1}$ has a fixed point.
Since $\lambda_{1}>1$ then  $\lim\limits_{w\to
  1}\omega(0,w)-\omega(0,F_{1}(w))>0.$
Proceeding as before we obtain:
\begin{center}$\liminf_{(x,y)\to (1,1)}
K_{\Delta^{2}}((0,0),(x,y))-\omega(0,f_{1}(x,y)) $
\\ $\leq\liminf_{(x,y)\to (1,1)}
K_{\Delta^{2}}((0,0),(F_{1}(w),w))-\omega(0,F_{1}(w))=\frac{1}{2}\log \lambda_{1} < +\infty$
\end{center} and by Julia's lemma for polydiscs $f_{1}$ admits restricted E-limit
$\tau_{1}=1.$
Then the inferior limit  ~\eqref{Julia per f_{1}} is bounded. Let  $\alpha$ denote
this inferior limit. This is a sort of  "boundary dilatation
coefficient" of $f_{1}$ in $(1,1).$ We can write:
\begin{center}$+\infty > \frac{1}{2}\log \alpha=\liminf\limits_{(x,y)\to(1,1)}
K_{\Delta^{2}}((0,0),(x,y))-\omega(0, f_{1}(x,y))=$
$=\frac{1}{2}\log \liminf\limits_{(x,y)\to
  (1,1)}\frac{1-|f_{1}(x,y)|}{1-|||(x,y)|||}= \frac{1}{2}\log
\liminf\limits_{t\to 1^{-}}\frac{1-|f_{1}(t,t)|}{1-t}.$
\end{center}
Let consider the holomorphic function  $\varphi:\Delta\to\Delta
$ defined by $\varphi(\xi)=f_{1}(\xi,\xi).$
By Lemma ~\ref{DilatCoeffHoros},
$\alpha_{\varphi}=\alpha_{f_{1}}.$
If $F_{1}(y_{0})=y_{0}$ it implies that
$\varphi(y_{0})=f_{1}(F_{1}(y_{0}),y_{0})=F_{1}(y_{0})=y_{0}$
then $\varphi(\xi)$ has a fixed point in $\Delta.$ We also know that
$\sigma(t)=(t,t)$ is a peculiar $(1,1)-$curve  then $f_{1}(t,t)\to 1$
when $t\to 1^{-}$. It follows that the point $1$ is a fixed point of
$\varphi$ on the boundary of $\Delta$ and since $\varphi$ has a fixed point then it must be
$\alpha_{f_{1}}=\alpha_{\varphi}>1$.
By Lemma ~\ref{DilatCoeffHoros}, this condition is sufficient to say that
the point $(1,1) $ cannot be a Wolff point for the map $f.$ In fact we
have:
$f_{1}(E((1,1),R))\subseteq E(1, \alpha_{f_{1}} R)\;\; \forall\;
R>0$
with $\alpha_{f_{1}}>1.$ It
implies
$f_{1}(E((1,1),R))\not\subseteq E(1,R)$
and the point $(1,1)$ cannot be a Wolff point of $f$.
\\
Suppose, now, $F_{1}$ has no fixed points. Then
$\liminf_{(x,y)\to(-1,-1)}
\max\{(\omega(0,x),\omega(0,y))\}-\omega(0,f_{1}(x,y))\leq
\liminf_{w\to -1}\omega(0,F_{1}(w))-\omega(0,F_{1}(w))=0,$
and then by  lemma ~\ref{Julia per f_{1}}, $f_{1}$ has restricted E-limit, say
$\tau\in\partial\Delta,$ in $
(-1,-1)$. Furthermore $\sigma_{1}(t)=(F_{1}(t),t)$ is a peculiar $(-1,-1)-$curve
and then, proceeding as before we have that $\tau=-1$ and
$f_{1}(E(-1,R)\times E(-1,R))\subseteq E(-1,R)\;\; \forall \;\; R>0.$
If the point $(1,1)$ is a Wolff point of $f$ we  have:
$f_{1}(E(1,R)\times E(1,R))\subseteq E(1,R)\;\; \forall \;\; R>0,$
and then, chosen $R>0$ such that $\{0\}\in \partial E(1,R) \cap
\partial E(1,R) $ we will have  $(0,0)\in \partial E((1,1),R))\cap \partial
 E((-1,-1),R))$
and thus $f_{1}(0,0)\in E(1,R)\cap E(-1,R)=\{0\}\Rightarrow f_{1}(0,0)=0
\Rightarrow F_{1}(0)=0.$
But it isn't possible because we supposed that  $F_{1}$ has
no fixed points in $\Delta$.
So, also in case $i)$ the point $(1,1)$ cannot be a Wolff point of $f$.
We can note, now, that the flat component of the boundary cannot be
a Wolff component of $f$ for otherwise $f_{1}(\cdot, y)$ or $f_{2}(x,\cdot)$
would have Wolff points but it isn't possible because, in this case, they have, both,
fixed points in $\Delta$.
By Lemma ~\ref{W(f)Cpa} we know that two points on the Silov boundary
of the same flat component cannot be, at the same time, Wolff point of
the map $f$. Furthermore $W(f)$ is arcwise connected, so we have just
one more possibility:
there is one point (different from $(1,1)$) on the Silov boundary
that is  the Wolff point of $f$.
But:

-  the point $(1,-1)$ cannot be Wolff points of $f$ because
  otherwise we would have
  $f_{2}(E(1,R) \times E(1,R))\subseteq  E(1,R)
 \;\; \forall\; R>0$ since $\lambda_{2}\leq 1$ and
$f_{2}(E(1,R) \times E(-1,R))\subseteq  E(-1,R)
 \;\; \forall\; R>0.$

Chosen $0\in \partial E(-1,R)\cap E(1,R)$ and $x\in E(1,R)$ then
$(x,0)\in [E((1,1),R)]\cap [E((-1,1),R)]$ and  we  obtain
$f_{2}(x,0)=0 \;\; \forall \; x\in E(1,R)$ and then
$F_{2}(x)\equiv 0$
but it is inconsistent with our hypothesis. So  the point   $(1,-1)$
cannot be Wolff points of $f.$

- the point $(-1,-1)$ cannot be Wolff point of $f$ because
  otherwise
$f_{2}(E(-1,R) \times E(-1,R))\subseteq  E(-1,R)
 \;\; \forall\; R>0.$
As above, chosen $0\in \partial E(-1,R)\cap E(1,R)$ we have
$f_{2}(0,0)=0$ that implies $ F_{2}(0)=0$ but it isn't possible because
$F_{2}$ has Wolff point $1$.

- the point  $(-1,1)$
cannot be Wolff point of $f$. This last statement follows by the point \emph{i)} of Proposition~\ref{R1R2}
Indeed if $(-1,1)$ were a Wolff point, we would have:
$f(E(-1,R)\times E(1,R))\subseteq E(-1,R)\times E(1,R)\;\; \forall \;
R>0.$
In particular, chosen $R_{1},R_{2}$ such that $\frac{\lambda_{12}}{R_{1}}=\frac{\lambda_{2}}{R_{2}}$,   $\exists\; x_{0}\in \partial E(-1,R)\cap \partial E(1,R_{1})$
such that $f(x_{0},y)=(x_{0},\tilde{y})$ for all $ y \in E(1,R_{2})$
 and
also such that $f_{1}(x_{0},y)=\nolinebreak[4]x_{0}$ $\forall \; y\in
E(1,R_{2}).$ It implies that $ F_{1}(y)=x_{0}$ $\forall \; y\in
E(1,R_{2})$. But it isn't possible.
So we proved that there isn't any Wolff point for $f$ and
in this way we end the proof of the point $i)$ of the theorem.

$iii), iv).$
If $f$  is of \emph{second type} it is clear that every point of the
flat component $\Gamma_{(1,0)}=\{1\}\times \Delta$ is a Wolff point, in
fact, chosen
$(1,\tilde{y})\in \Gamma_{(1,0)},$ the small horosphere centered in this point  is
$E((1,\tilde{y}),R)=E(1,R)\times \Delta\;\; \forall\; \tilde{y}\in
\Delta$ and we have:
\begin{center}
$f(E((1,\tilde{y}),R))=(f_{1}(E(1,R)\times \Delta),f_{2}(E(1,R)\times
\Delta))$
\end{center}
but $f_{1}(E(1,R)\times \Delta)\subseteq E(1,R)$ because
  $\tau=1$ is the Wolff point of $f_{1}(\cdot,y)\;\; \forall\; y\in\Delta$
and $f_{2}(E(1,R)\times\Delta)\subseteq \Delta.$ It follows that
$$f(E((1,\tilde{y}),R))\subseteq E(1,R)\times \Delta=E((1,\tilde{y}),R)\;\; \forall\; \tilde{y}\in
\Delta$$
and then we have that every point of the flat component $\Gamma_{(1,0)} $
is a Wolff point of $f.$
Furthermore, applying again Julia's lemma for polydiscs, as in the proof of point $ii)$ of this 
theorem, we have that  $\lambda_{2}\leq 1$ implies 
$f_{1}(E((1,1),R))\subseteq E((1,1),R)\;\; \forall\; R>0
\;\; \hbox{and}\;\; \{(1,1)\}\in W(f).$
If on the other hand  $\lambda_{2}>1,$ by the proof of the points
\emph{$i)$} and \emph{ $ ii)$} of the theorem,  we have that $\{(1,1)\}\notin W(f).$
Also in these cases we need to prove that there isn't any other Wolff point.
If $\lambda_{2}<1$
the points of the flat component  $\Delta\times\{1\}$ cannot be  Wolff
points for $f$ because $f_{2}$ has
fixed points.
Moreover by lemma ~\ref{W(f)Cpa}, the points of the Silov boundary of $\Gamma_{01}$
cannot be Wolff point because $\{(1,1)\}\in W(f).$ Then since $W(f)$ is arcwise connected there isn't any other possibility
and the set of the Wolff points of $f$ must be $W(f)=\{\{1\} \times \Delta\}\cup \{(1,1)\}$.
If $\lambda_{2}>1,$
  neither the points of the flat component  $\Gamma_{(0,1)}$
nor the points of  $\Gamma_{(0,-1)}$ can be Wolff points because $f_{2}$
has fixed points.
Furthermore we have $\{(1,-1)\}\notin W(f).$ Indeed
$f_{1}(E(1,R)\times\Delta)\subseteq E(1,R)\;\; \forall\; R>0 $ and
$f_{2}(E(1,R)\times E(1,R)) \subseteq E(1, R)\;\; \forall\; R>0.$
Thus, using Proposition ~\ref{R1R2}, we can
prove that  $\{(1,-1)\}\notin W(f).$

\emph{$v)$} As in \emph{$iii)$} and \emph{$iv)$}, we have that every point
of the flat components of the boundary, $\{1\}\times\Delta$ and
$\Delta\times\{1\},$ is a
Wolff point of $f$. It implies
$\{(1,1)\}\in W(f),$ indeed
$f(E((1,1),R))=\big{(}f_{1}([E(1,R)\times\Delta]\cap[\Delta\times
E(1,R)]),f_{2}([E(1,R)\times\Delta]\cap[\Delta\times
E(1,R)]) \big{)}\subseteq$
$\subseteq E(1,R)\times E(1,R)=E((1,1),R).$
Thus $W(f)\neq \emptyset $ and it is at least
$[\{1\}\times\Delta]\cup\{(1,1)\}\cup [\Delta\times\{1\}].$
By lemma ~\ref{W(f)Cpa}, it follows that it must be exactly $ W(f)=[\{1\}\times\Delta]\cup\{(1,1)\}\cup [\Delta\times\{1\}] .$
In this way we proved every statement of the theorem in the direction
 $\Leftarrow$. Then we can conclude that also the implications in direction
$\Rightarrow$ are proved.

\emph{$vi)$}
Let suppose that $f_{1}(x,y)=x\;\; \forall\; y\in \Delta$ and let
examine the set $W(f)$.
Let note that $f_{2}\neq\pi_{2}$ otherwise $f=id_{\Delta}.$
Furthermore there is no a holomorphic map $F_{2}:\Delta\to \Delta $ such that
  $f_{2}(x,F_{2}(x))=F_{2}(x)\;\; \;\forall\; x\in
  \nolinebreak[4]\Delta.$ Otherwise we will have that 
$\exists\; x_{0}\in\Delta$ such that $f(x_{0},F_{2}(x_{0}))=(f_{1}(x_{0},F_{2}(x_{0})),f_{2}(x_{0},F_{2}(x_{0})))=(x_{0},F_{2}(x_{0}))$
and it is not possible because we supposed $f$ without fixed points in
the bidisc.
Then, by Herv{\'e} theorem,  $f_{2}$ has Wolff point $\tau_{f_{2}}=e^{\imath\theta_{2}}.$  Let suppose, without loss of generality, that
$e^{\imath\theta_{2}}=1.$ Then every point of
$\Gamma_{(0,1)}$ is a Wolff point of $f$, indeed, by $f_{1}=id$ and $\tau_{f_{2}}=1$
we obtain:
$f(\Delta\times
E(1,R))=(f_{1}(\Delta\times
E(1,R)),f_{2}(\Delta\times
E(1,R)))\subseteq \Delta\times E(1,R)\;\forall\;R>0$
and, also:
$f(E(1,R)\times\Delta)=(f_{1}(E(1,R)\times\Delta),f_{2}(E(1,R)\times\Delta))\subseteq E(1,R)\times\Delta\;\forall\;R>0.$
It follows immediately that every point of
$\Gamma_{(1,0)}$ is a Wolff point and then, 
the point $\{(1,1)\}$ is a Wolff point of $f$.
In the same way we can prove that $\Gamma_{(-1,0)}$
is a Wolff component of $f$, 
and then also the point  $\{(-1,1)\}\in W(f).$
Thus $f$ fixes every leaf $(k,y)$:
$f(k,y)=(f_{1}(k,y),f_{2}(k,y))=(k,\tilde{y})$ and in particular $y\in E(1,R)$ implies $ \tilde{y}\in E(1,R).$
By this remark follows that the points of the flat component of the boundary
$\Gamma_{(0,-1)}=\Delta\times\{-1\}$ cannot be  Wolff points for $f$
otherwise we will have:
$f_{2}(\Delta\times E(1,R)) \subseteq  E(1,R)\;\; \forall\; R>0$ and
$f_{2}(\Delta\times E(-1,R)) \subseteq  E(-1,R)\;\; \forall\; R>0.$
Then chosen $\{0\}\in \partial E(1,R)\cap \partial E(-1,R)$ we have
$f_{2}(\Delta\times E(1,R))\subseteq [E(1,R)\cap E(-1,R)]=\{0\}\;\;\hbox{and} ;
f_{2}(x,0)=x\;\; \forall\; x\in \Delta.$
Consequently:
$f(0,0)=(f_{1}(0,0),f_{2}(0,0))=(0,0)$ but it is inconsistent with
the hypothesis that $f$ has not fixed point in $\Delta^{2}$.
Also the points  $\{(-1,-1)\}$ e $\{(1,-1)\}$
cannot be Wolff points of $f,$ because we already know that
$\{(1,1)\}$ and $\{(-1,1)\}$ are Wolff points of $f$.
In fact if $\{(-1,-1)\}\in W(f)$ were  we would have:
$f_{2}(E(-1,R)\times E(-1,R)) \subseteq E(-1,R)\;\;\forall\; R>0$ and 
 $f_{2}(E(1,R)\times E(1,R)) \subseteq E(1,R)\;\;\forall\; R>0. $
As done before, we can choose $\{0\}\in \partial E(1,R)\cap \partial
E(-1,R)$ and we obtain:
$f(0,0)=(f_{1}(0,0),f_{2}(0,0))=(0,0)$ but it isn't possible.
In the same way we can prove that $\{(-1,1)\}\notin W(f).$
Then:
$W(f)=\{\{-1\} \times \Delta \}\cup \{(-1,1)\} \cup \{\Delta
\times \{1\}\}\cup \{(1,1)\} \cup \{\{1\} \times \Delta \}$
and it ends the proof of the theorem.
\end{proof}
We can end this section proving theorem ~\ref{FIX1}. Recall that by hypothesis $f_{i}(0,0)=(0,0)$ and $F_{i}(0)=0$, $i=1,2.$
Moreover if we denote by $K((0,0),R)$ the kobayashi disk centered in $(0,0)$ with radius $R,$ we have that $f(K((0,0),R))\subseteq K((0,0),R)$ for all $R>0.$

\begin{proof}\emph{Theorem ~\ref{FIX1}. }
\begin{remark}
Let denote by $\Gamma_{(x,y)}$ the flat component of
$\partial\Delta^{2},$ containing the point $(x,y).$ \end{remark}
 Recall that by hypothesis $f_{i}(0,0)=(0,0)$ and $F_{i}(0)=0$, $i=1,2.$
Moreover if we denote by $K((0,0),R)$ the Kobayashi disk centered
in $(0,0)$ with radius $R,$ we have that $f(K((0,0),R))\subseteq
K((0,0),R)$ for all $R>0.$

Suppose first that $dm Fix(f)=0.$
If a point of the Silov
boundary, say $(1,1),$ is a Wolff point then $f(E(1,1),R)\subseteq
(E(1,1), R)$ for all $R>0.$ If we take $(0,0)\neq (x_{0},y_{0})=
[\partial(E(1,1), R)\cap
\partial K((0,0),R_{1})]$ we have that   $(x_{0},y_{0})$ must be
fixed by $f$ and it is a contradiction. On the other hand if a
point of a flat component, say $ (0,1)$, is a Wolff point, every
point of that flat component is a Wolff point for $f.$ In this
case we have $f(\Delta \times E(1,R) )\subseteq \Delta\times
E(1,R) $ for all $R>0.$ Chosen $(x, y_{0})\in [\partial
\Delta\times E(1,R)]\cap
\partial K((0,0),R)$ we obtain that $f(x, y_{0})=(x_{1}, y_{0})$
for every $x\in \Delta.$ Then $f_{2}(x,y)=y$but it a contradiction
because, in this case $\dim Fix(f)=0$ and also $\dim Fix f_{2}=0.$
Thus we conclude that $W(f)=\emptyset.$

Now, suppose that $dim Fix{f}=1.$

1) If $(e^\imath\theta,1)\in W(f)$ then $f(E((e^\imath\theta ,1),
R))\subseteq E((e^\imath\theta,1),R)$ for all $R>0.$  Since, by
definition, $K((0,0), R_{1})$ is the product of two Poincar{\'e}
discs of radius $R_{1},$ we can take $(0,0)\neq (e^\imath\theta
x_{0},x_{0})= [\partial(E(e^\imath\theta,1), R)\cap \partial
K((0,0),R_{1})).$ We have that $(e^\imath\theta x_{0},x_{0})$ must
be fixed by $f.$ We can do the same thing for every point
$(e^\imath\theta x, x)\in \Delta^{2},$ changing the radius
$R_{1}.$ Then $Fix(f)$ is equal to the geodesic $(e^\imath\theta
z,z).$ Using the estimate of the inferior limit in lemma
~\ref{JuliaForPolydiscs} as in the points of theorem ~ref{W(f)} we
see that $\partial G =W(f).$ Indeed any other point of the flat
component of $\partial Delta^{2}$ cannot be a Wolff point since
$f_{i}, i=1,2$ has fixed points. On the other hand if $g(z)=
e^\imath\theta z$ then $(e^\imath\theta,1)\in W(f)$ and every
point of the boundary of $\partial G \in W(f).$ To prove this fact
it is sufficient to apply lemma ~\ref{JuliaForPolydiscs} and study
the inferior limit ~\ref{limiteJulia} along the direction
$(e^\imath\theta t,t)$, as in the proof of the theorem
~\ref{W(f)}.

2), 3) Suppose that $g$ is proper and $g$ is neither an
automorphism nor the identity. Then $W(f) $ is contained in the
Silov boundary of $\Delta^{2}$ and it is in contradiction with the
preceding point. Then $W(f)=\emptyset.$ It proves the implication
"$\Rightarrow$ " of point $2.$ Suppose now that $g$ is not a
proper map and let prove that $W(f)\neq\emptyset$ and it is
disconnected. Since $g$ is not proper then there exists a sequence
$z_{k}\in \Delta$ such that $z_{k}\to
e^{\imath\theta}\in\partial\Delta $ as $k\to\infty$ and
$g(z_{k})\to c\in\Delta$ as $k\to\infty.$ By theorem
~\ref{JuliaForPolydiscs} we get that
\[
\liminf\limits_{(x,y)\to(e^{\imath\theta},c)}
K_{\Delta^{2}}((0,0),(x,y))-\omega(0, f_{i}(x,y))\leq\]
\[
\leq\liminf\limits_{k\to\infty}
K_{\Delta^{2}}((0,0),(g(z_{k}),z_{k}))-\omega(0,
f_{i}(g(z_{k}),z_{k}))=\] \[ \liminf\limits_{k\to\infty} \max
\{\omega(0,g(z_{k}));\omega(0,z_{k})\}-\omega(0,
f_{i}(g(z_{k}),z_{k}))=(\star)\]
if $i=1$ then $(\star)=\infty;$ and  if $i=2$ then $(\star)=0.$

 Thus we get
 \[f(E(e^{\imath\theta},c),R)\subseteq (E(e^{\imath\theta},c),R)\;
 \forall\; R>0 \]
 and $(e^{\imath\theta},c)\in W(f)\neq \emptyset.$
 Thus $\Gamma_(e^{\imath\theta},c)\in W(f).$ It means that $f_{2}(x,\cdot)$
 has Wolff point but  since $f_{2}$ has also fixed points it follows that $f_{2}=\pi_{2}.$
It follows also that the flat component of the boundary $\Gamma_(e^{-\imath\theta},c)\in W(f).$
Any other point of the flat component can be a Wolff point since $f_{1}\neq\pi_{1}$ and $f_{1}(\cdot,y)$
has not Wolff points.
Any other point of the Silov boundary can be a Wolff point since the point $1)$
of the theorem holds. Thus $W(f)$ is disconnected.
In this way we proved implication $\Rightarrow$ of point $3).$
Now we have that if $W(f)=\emptyset$ then $g$ is proper for the preceding point.
And if $f_{2}=\pi_{2}$ then $g$ is not proper since $W(f)\neq id.$
\end{proof}

%
%
%
%
\section{Examples}
We can now give an example for each case of the theorem:
\\
\emph{Example $i)$ of theorem~\ref{W(f)}}
We are going to give an example of a holomorphic self map $f,$ of the
complex bidisc, without fixed points and without Wolff points.
Let consider:
$f(x,y)=(\frac{1}{2}(x+F_{1}(y)),\frac{1}{2}(y+F_{2}(x)))$
with  $F_{1}$ and $F_{2}$ holomorphic self map of the unit disc
$\Delta.$
We can choose
$F_{1}$ and $F_{2}$ such that
$f$ hasn't fixed points in  $\Delta^{2}$ that is for example:
$F_{1}(y)=y^{2}$ and $F_{2}(x)=\frac{3x+1}{x+3}.$
We are going to prove that the point $(1,1)$ isn't Wolff point for $f.$
By lemma ~\ref{DilatCoeffHoros}  it will be sufficient prove that the limit
(~\ref{limiteJulia}) is strictly greater than 0. 
If $\max\{(\omega(0,x),\omega(0,y))\}=\omega(0,x):$
\begin{center}$ \liminf\limits_{(x,y)\to(1,1)}
\max\{(\omega(0,x),\omega(0,y))\}-\omega(0,f_{1}(x,y))=\liminf\limits_{(x,y)\to(1,1)}
\omega(0,x)-\omega(0,\frac{1}{2}(x+F_{1}(y)))$ \\$\geq
\liminf\limits_{(x,y)\to(1,1)}  \frac{1}{2}\log\Big[
\frac{1+|x|}{1+|\frac{1}{2}(x+F_{1}(y))|}\;\;\frac{1-|\frac{1}{2}x|-|\frac{1}{2}F_{1}(y)|}{1-|x|}\Big]
>\liminf\limits_{(x,y)\to(1,1)}  \frac{1}{2}\log\Big[
\frac{1+|x|}{1+|\frac{1}{2}(x+F_{1}(y))|}\;\;\;\frac{1-|\frac{1}{2}x|-|\frac{1}{2}x|}{1-|x|}\Big]=0.$
\end{center}
Then the boundary dilatation coefficient of $f_{1}$
at the point 
$ z=1$ is  $\alpha_{f_{1}}>1.$
On the other hand if $\max\{(\omega(0,x),\omega(0,y))\}=\omega(0,y)$
we have:
\begin{center}
$\liminf\limits_{(x,y)\to(1,1)}
\max\{(\omega(0,x),\omega(0,y))\}-\omega(0,f_{1}(x,y))=\liminf\limits_{(x,y)\to(1,1)}
\omega(0,y)-\omega(0,\frac{1}{2}(x+F_{1}(y)))$
$\geq
\liminf\limits_{(x,y)\to(1,1)}\frac{1}{2}\log \Big[
\frac{1+|y|}{1+|\frac{1}{2}(x+F_{1}(y))|}\;\;\;\frac{1-|\frac{1}{2}x|-|\frac{1}{2}F_{1}(y)|}{1-|y|}\Big]$
$>\liminf\limits_{(x,y)\to(1,1)}\frac{1}{2}\log \Big[
\frac{1+|y|}{1+|\frac{1}{2}(x+F_{1}(y))|}\;\;\;\frac{1-|\frac{1}{2}y|-|\frac{1}{2}y|}{1-|y|}\Big]=0.$
\end{center}
Then  $\forall \;\;R>0 \;\;f_{1}(E((1,1),R))\subseteq E(1,\alpha_{f_{1}} R)$
with $\alpha_{f_{1}}>1$ and consequently, by the previous proof
$(1,1)$ cannot be a Wolff point for $f.$ Furthermore, by the proof of
$i)$ we can conclude that there isn't any other Wolff point.\\
\emph{Example $ii)$ of theorem~\ref{W(f)}}: 
$f(x,y)=(\frac{1}{2}(x+\frac{5y+3}{3y+5}), \frac{1}{2}(y+\frac{3x+1}{x+3})).$
\\
\emph{Example $iii)$ of theorem~\ref{W(f)}}:
$f(x,y)=(\frac{3x+1}{x+3}, \frac{1}{2}(y+F_{2}(x)))$ with
$F_{2}(x)=\frac{5x+3}{3x+5}.$
\\
\emph{Example $iv)$ of theorem~\ref{W(f)}}:
$f(x,y)=(\frac{3x+1}{x+3}, \frac{1}{2}(y+F_{2}(x)))$ with
$F_{2}(x)=x^{2}.$
\\
\emph{Example $v)$ of theorem~\ref{W(f)}}
$f(x,y)=(\frac{3x+1}{x+3},\frac{5y+3}{3y+5} ).$

\end{document}